\documentclass[12pt,reqno,a4paper,twoside]{article}
\usepackage{amsmath,amsthm,amstext,amscd,amssymb,euscript}
\usepackage{epsf,color}
\renewcommand{\phi}{\varphi}

\newcommand{\Zd}{\mathbb Z^d}

\def\1{ {\mathit{1} \!\!\>\!\! I} }

\makeatletter
\@addtoreset{equation}{section}
\makeatother

\newtheorem{ittheorem}{Theorem}
\newtheorem{itlemma}{Lemma}
\newtheorem{itcorollary}{Corollary}
\newtheorem{itproposition}{Proposition}
\newtheorem{itdefinition}{Definition}
\newtheorem{itremark}{Remark}

\newenvironment{theorem}{\addtocounter{equation}{1}
\begin{ittheorem}}{\end{ittheorem}}

\newenvironment{lemma}{\addtocounter{equation}{1}
\begin{itlemma}}{\end{itlemma}}
\newenvironment{corollary}{\addtocounter{equation}{1}
\begin{itcorollary}}{\end{itcorollary}}

\newenvironment{proposition}{\addtocounter{equation}{1}
\begin{itproposition}}{\end{itproposition}}

\newenvironment{definition}{\addtocounter{equation}{1}
\begin{itdefinition}}{\end{itdefinition}}

\newenvironment{remark}{\addtocounter{equation}{1}
\begin{itremark}}{\end{itremark}}

\newcommand{\beq}{\begin{eqnarray}}
\newcommand{\eeq}{\end{eqnarray}}

\newcommand{\be}{\begin{equation}}
\newcommand{\ee}{\end{equation}}

\newcommand{\bl}{\begin{lemma}}
\newcommand{\el}{\end{lemma}}

\newcommand{\br}{\begin{remark}}
\newcommand{\er}{\end{remark}}

\newcommand{\bt}{\begin{theorem}}
\newcommand{\et}{\end{theorem}}

\newcommand{\bd}{\begin{definition}}
\newcommand{\ed}{\end{definition}}

\newcommand{\bp}{\begin{proposition}}
\newcommand{\ep}{\end{proposition}}

\newcommand{\bc}{\begin{corollary}}
\newcommand{\ec}{\end{corollary}}

\newcommand{\bpr}{\begin{proof}}
\newcommand{\epr}{\end{proof}}

\newcommand{\bi}{\begin{itemize}}
\newcommand{\ei}{\end{itemize}}

\newcommand{\ben}{\begin{enumerate}}
\newcommand{\een}{\end{enumerate}}

\newcommand{\R}{\mathbb R}
\newcommand{\N}{\mathbb N}

\newcommand{\E}{\mathbb E}

\newcommand{\pee}{\ensuremath{\mathbb{P}}}

\newcommand{\loc}{\ensuremath{\mathcal{L}}}

\newcommand{\la}{\ensuremath{\Lambda}}

\newcommand{\epsi}{\ensuremath{\epsilon}}

\newcommand{\sip}{SIP(m)}

\def\now{
\ifnum\time<60
          12:\ifnum\time<10 0\fi\number\time am
          \else
            \ifnum\time>719\chardef\a=`p\else\chardef\a=`a\fi
          \hour=\time
          \minute=\time
          \divide\hour by 60 
          \ifnum\hour>12\advance\hour by -12\advance\minute by-720 \fi
          \number\hour:%
          \multiply\hour by 60 
          \advance\minute by -\hour
          \ifnum\minute<10 0\fi\number\minute\a m\fi}
\newcount\hour
\newcount\minute
\numberwithin{equation}{section}         

\theoremstyle{remark}

\newcommand{\caD}{{\mathcal D}}

\newcommand{\caP}{{\mathcal P}}

\newcommand{\laa}{\overline{\lambda}}
\newcommand{\six}{\sum_{i=1}^n \delta_{x_i}}
\newcommand{\siP}{SIP(1)}
\newcommand{\roo}{\overline{\rho}}

\begin{document}
\title{{\bf Correlation inequalities for interacting particle
systems with duality}}
\author{
C. Giardin\`a
\footnote{Modena and Reggio Emilia University, viale Allegri 9, 42121
Reggio Emilia, Italy, \newline {\em cristian.giardina@unimore.it}}\\
F.\ Redig
\footnote{Universiteit Nijmegen, IMAPP,
Heyendaalseweg 135
6525 AJ Nijmegen
The Netherlands,
{\em redig@math.leidenuniv.nl}}\\
K.\ Vafayi\footnote{Mathematisch Instituut Universiteit Leiden,
Niels Bohrweg 1, 2333 CA Leiden, The Netherlands,
{\em vafayi@math.leidenuniv.nl}}
}
\maketitle
\begin{quote}
Abstract: We prove a {\em comparison inequality} between a system of
independent random walkers and a system of random walkers which
either interact by {\em attracting each other} -- a process which we
call here the symmetric inclusion process (SIP) -- or {\em repel each
other} -- a generalized version of the well-known symmetric exclusion
process. As an application, new {\em correlation inequalities} are
obtained for the SIP, as well as for some interacting diffusions
which are used as models of heat conduction, -- the so-called Brownian
momentum process, and the Brownian energy process. These
inequalities are counterparts of the inequalities (in the opposite
direction) for the symmetric exclusion process, showing that the SIP
is a natural bosonic analogue of the symmetric exclusion process,
which is fermionic. Finally, we consider a boundary driven version
of the SIP for which we prove duality and then obtain correlation
inequalities.
\end{quote}
\vspace{12pt}

\section{Introduction}
In Liggett \cite{ligg},
Chapter VIII, proposition 1.7, a comparison inequality between
independent symmetric random walkers and corresponding exclusion
symmetric random walkers is obtained. This inequality plays a crucial role in
the understanding of the exclusion process (SEP); it makes rigorous the
intuitive picture that symmetric random walkers interacting by
exclusion are more spread out than the corresponding independent
walkers, as a consequence of their repulsive interaction (exclusion),
or in more physical terms, because of the fermionic nature of the
exclusion process. The comparison inequality is a key ingredient in
the ergodic theory of the symmetric exclusion process, i.e., in the
characterization of the invariant measures, and the
measures which are in the course of time attracted to a given
invariant measure. The comparison inequality has been generalized
later on by Andjel \cite{Andjel}, Liggett \cite{ligg2},  and
recently in the work of  Borcea, Br\"{a}nd\'{e}n and Liggett \cite{brand}.

In the search of a natural conservative particle system where the
opposite inequality holds, i.e., where the particles are {\em less
spread out} than corresponding independent random walkers, it is
natural to think of a ``bosonic counterpart'' of the exclusion
process. In fact, such a process was introduced in
\cite{gkr} and \cite{gkrv} as the {\em dual} of the
Brownian momentum process, a stochastic model of heat conduction
(similar models of heat conduction were introduced in \cite{bol} and
\cite{gk}, see also \cite{ber} for the study of the structure
function in a natural asymmetric version).

In the present paper we analyze this ``bosonic counterpart'' of the
exclusion process. We will call this process (as will be motivated by
a Poisson clock representation) the ``symmetric inclusion process'' (SIP).
In the SIP, jumps are performed according to independent random walks, and
on top of that particles ``invite'' other particles to join their
site (inclusion). For this process we prove the
analogue of the comparison inequality for the symmetric exclusion
process. From the comparison inequality, using the knowledge of the
stationary measure and the self-duality property of the process, we
deduce a series of correlation inequalities. Again, in going from
exclusion to inclusion process the correlations turn from negative
to positive. We remark however that these positive
correlation inequalities are different from the ordinary
preservation of positive correlations for monotone processes
\cite{harris}, because the SIP is not a monotone process.
Since the SIP is dual to the heat conduction model
it is immediate to extend those correlation inequalities to the
Brownian momentum process and the Brownian energy process.

We also introduce the non-equilibrium versions of
the SIP, i.e., we consider the boundary driven version of SIP. In this case, for
appropriate choice of the boundary generators, we
prove duality of the process to a SIP model with absorbing boundary
condition. We then deduce a correlation inequality, explaining and
generalizing the positivity of the covariance in the non-equilibrium
steady state of the heat conduction model in \cite{gkr}.

All the results will be stated in the context
of a family of SIP$(m)$ models, which are labeled by parameter $m\in\N$.
As the SEP model can be generalized to the situation where there are
at most $n\in\N$ particles per site (this corresponds to a quantum
spin chain with SU(2) symmetry and spin value $j=n/2$), in the same
way the SIP model can be extended to represent the situation of a
quantum spin chain with SU(1,1) symmetry and spin value $k=m/4$ \cite{gkr}.

The paper is organized as follows. In Section 2 we define the
SIP$(m)$ process, restricting to a context where its existence can
be immediately established. The main comparison inequality, which
allows to compare SIP walkers to independent walkers (by a suitable
generalization of Liggett comparison inequality) is proved in
Section 3. Correlation inequalities for the SIP$(m)$ process that
can be deduced from the comparison inequality are proved in Section
5 (the necessary knowledge of the stationary measure and the
self-duality property are presented in Section 4). In particular, in
Section 5 it is proved that when the SIP$(m)$ process is started
from its stationary measure then correlations are always positive,
while when the process is initialized with a general product measure
then positivity of correlations is recovered in the long time limit.
Further correlation inequalities for systems similar to the SIP$(m)$
process are discussed in the subsequent Sections. Attractive
interaction (the SEP$(n)$, which generalize the standard SEP) is
presented in Section 6. Some interacting diffusions dual to the
SIP$(m)$ process are studied in Section 7. Finally the boundary
driven SIP$(m)$ process is analyzed in Section 8.

\section{Definition}
In the whole of the paper, $S$ will denote or a finite set, or
$S=\Zd$. Next, $p(x,y)$ denotes an irreducible (discrete-time)
symmetric random walk transition probability on $S$, i.e.,
$p(x,y)=p(y,x)\geq 0$, $\sum_y p(x,y)= 1$, and $p(x,x)=0$. In the
case  $S=\Zd$,  we suppose furthermore that $p(x,y)$ is finite range
and translation invariant, i.e., $p(x,y) =\pi(y-x)$, and there
exists $R>0$ such that $p(x,y)=0$ for $|x-y|>R$. This assumption for
the infinite-volume case avoids technical problems for the existence
of the $\sip$ which for the subject of this paper are irrelevant.
The proof of existence of the $\sip$ in our infinite-volume context
(with the process started from a ``tempered'' initial configuration, i.e. $\eta(y) \le ||y||^k$ for some
$k$ and for all $y$) follows from self-duality, along the lines of \cite{demasi}, Chapter
2.

The symmetric inclusion process with parameter $m\in (0,\infty)$ associated to the transition kernel $p$ is the
Markov process
on $\Omega:=\N^S$ with generator
defined on the core of local functions by
\be\label{sipgen1}
L f(\eta) = \sum_{x,y\in S} p(x,y) 2\eta_x (m+2\eta_y) \left(f(\eta^{x,y}) -f(\eta)\right)
\ee
where, for $\eta\in\Omega$, $\eta^{x,y}$ denotes the configuration obtained from $\eta$ by removing
one particle from $x$ and putting it at $y$.

In \cite{gkr}, for $m=1$ this model was introduced as the dual of a model of
heat conduction, the so-called Brownian momentum process, see
also \cite{gkrv}, and \cite{bol} for generalized and or similar models of
heat conduction.

The process with generator \eqref{sipgen1} can be interpreted as
follows. Every particle has two exponential clocks: one clock -the
so-called random walk clock- has rate $2m$, the other clock -the
so-called inclusion clock- has rate 4. When the random walk clock of
a particle at site $x\in S$ rings, the particle performs a random
walk jump with probability $p(x,y)$ to site $y\in S$. When the
inclusion process clock rings at site $y\in S$, with probability
$p(y,x)=p(x,y)$ a particle from site $x\in S$ is selected and joins
site $y$.

From this interpretation, we see that besides jumps of a system of independent
random walkers, this system of particles has the tendency to
bring particles together at the same site (inclusion), and can therefore
be thought of as a ``bosonic'' counterpart of the symmetric exclusion
process.

To make the analogy with the exclusion process even more
transparent, in an exclusion process with at most
$n$ particles ($n\in\N$) per site (notation $SEP(n)$), the jump rate is
$\eta_i(n-\eta_j) p(i,j)$. Apart from a global factor 4, the $SIP(m)$ is obtained by changing the minus
into a plus and choosing $n=m/2$.

Notice that the rates in \eqref{sipgen1}
are increasing both in the number of
particles of the departure and in
the number of
particles of the
arrival site
(the rate is $p(x,y) 2\eta_x(m+2\eta_y)$ for
a particle to jump from $x$ to $y$).
Therefore, by the necessary and sufficient
conditions of \cite{gs}, Theorem 2.21, the SIP
is not a monotone process. It is also easy to see that due to the
attraction between particles in the SIP, there cannot be
a coupling that preserves the order of configurations, i.e.,
in any coupling starting from an unequal ordered pair of configurations,
the order will be lost in the course of time with
positive probability.

\subsection{Assumptions on the transition probability kernel}

In this section we introduce the assumptions that
we need to prove the positivity of correlations of stationary
measures obtained as limits of general initial product measures
(see later for precise definitions). This assumptions are only relevant
in the infinite volume case $S=\Zd$ and they are indeed satisfied in the
context of finite-range translation-invariant underlying random walk kernel
$p(x,y)= \pi (y-x)$.
However, all our results on correlation inequalities for
stationary measures depend only on one or both of
the assumptions below, i.e., if on more general
graphs, or on $\Zd$ with more general $p(x,y)$, existence
of $\sip$ would be established, then the corresponding correlation
inequalities hold under one or both of the assumptions A1, A2 below.

We define the associated continuous-time
random walk transition probabilities of random walk jumping at rate
$2m$:
\be\label{coco}
 p_{t} (x,y) = \sum_{n=0}^\infty \frac{(2mt)^n}{n!} e^{-2mt} p^{(n)} (x,y)
\ee
where $p^{(n)}$ denotes the $n^{th}$ power of
the transition matrix $p$. Denote by $\pee^{IRW(m)}_{x,y}$ the
probability measure on path space associated to two independent
random walkers $X_t,Y_t$ started at $x,y$ and jumping according to
\eqref{coco} and by $\pee^{SIP(m)}_{x,y}$ the corresponding probability
for two SIP walkers $X^{'}_t,Y^{'}_t$ jumping with the rates of
generator (\ref{sipgen1}).

We consider two assumptions
\begin{itemize}
\item[-] Assumption (A1) \be\label{rwcon} \lim_{t\to\infty}
\sup_{x,y}\pee^{IRW(m)}_{x,y} (X_t=Y_t)=0 \ee
\item[-] Assumption (A2) \be\label{rwcond} \lim_{t\to\infty}
\sup_{x,y}\pee^{SIP(m)}_{x,y} (X^{'}_t=Y^{'}_t)=0 \ee
\end{itemize}

The assumption (A1) amounts to requiring that for
large $t>0$, two independent random walkers walking according to the
continuous time random walk probability (\ref{coco}) will be at the
same place with vanishing probability. The assumption (A1) follows
immediately if we have
\be\label{rwcondo}
 \lim_{t\to\infty}\sup_{x,y} p_t(x,y)=0
\ee since then \be\label{rwconard}
 \lim_{t\to\infty}\sup_{x,y}\pee^{IRW(m)}_{x,y} \left( X_t= Y_t\right)
=\lim_{t\to\infty}\sup_{x,y}\sum_{u\in S} p_t
(x,u)p_t(y,u)=\lim_{t\to\infty}\sup_{x,y} p_{2t} (x,y)=0 \ee
Notice also that, by simple rescaling of time, (A1) holds for all $m>0$ as soon as it
holds for some $m>0$.

Assumption (A2) guarantees that
two walkers evolving with the SIP dynamic will be typically at
different positions at large times.
Notice that in the case we consider, i.e., the translation invariant finite-range case
$S=\Zd$, $p(x,y) = p(0,y-x)=: \pi (y-x)$, this is automatically
satisfied, as the difference walk $X^{'}_t-Y^{'}_t$ of two SIP
particles is a random walk $Z_t$ on $\Zd$ with generator
\be
 L^Z f(z) = 8\pi (z) (f(0)- f(z)) + \sum_{y} 4m\pi (y) (f(z+y)- f(z))
\ee
which is clearly not positive recurrent.

Assumption (A2) implies that any finite
number of SIP particles will eventually be at different locations.
This is made precise in Lemma \ref{vanishlem} in section
\ref{corre}.

\section{Comparison of the SIP with independent random walks}

We will first consider the SIP process with a finite number of
particle in subsection \ref{finite} and then state the comparison
inequality in subsection \ref{compare}.

\subsection{The finite SIP}\label{finite}

If we start the SIP with $n$ particles
at positions $x_1,\ldots,x_n\in S$, we can keep
track of the labels of the particles. This gives
then a continuous-time Markov chain
on $S^n$ with generator
\beq\label{finsipgen}
\loc_n f(x_1,\ldots,x_n) &=& \sum_{i=1}^n\sum_{y\in S}
2p(x_i,y)\left(m+2\sum_{j=1}^n I(y=x_j)\right) (f(x^{x_i,y})-f(x))
\nonumber\\
&=&
\loc_{1, n} f (x) + \loc_{2,n} f (x)
\eeq
where $x^{x_i,y}$ denotes
the $n$-tuple $(x_1,\ldots,x_{i-1},y,x_{i+1},\ldots,x_n)$. Further,
$\loc_{1,n}$, resp.\ $\loc_{2,n}$ denote the
random walk resp.\  inclusion part of the generator and are
defined as follows
\be\label{walkpart}
\loc_{1, n} f(x_1,\ldots,x_n) = \sum_{i=1}^n\sum_{y\in S}
2mp(x_i,y)(f(x^{x_i,y})-f(x))
\ee
\be\label{inpart}
\loc_{2,n} f(x_1,\ldots,x_n) = \sum_{i=1}^n\sum_{j=1}^n
4p(x_i,x_j)(f(x^{x_i,x_j})-f(x))
\ee

\subsection{Comparison inequality}\label{compare}

From the description above, it is intuitively clear that in the SIP,
particle tend to be less spread out than in a system of independent
random walkers. Theorem \ref{siplig} below formalizes this
intuition and is the analogue of a comparison inequality of the SEP
(\cite{ligg}, Chapter VIII, Proposition 1.7).

To formulate it, we need the notion of a positive definite function.
A function $f: S\times S\to \R$ is called positive definite if for
all $\beta:S\to\R$ such that $\sum_{x} |\beta(x)| < \infty$
\[
 \sum_{x,y} f(x,y)\beta(x)\beta(y) \geq 0
\]
A function $f:S^n\to\R$ is called positive definite if it is
positive definite in every pair of variables.

We first introduce a slightly more general generator with parameters
$a>0$, $b\in \R$ that includes
both process of exclusion and inclusion type.
\be\label{abgen}
\loc^{a,b}_n f(x_1,\ldots,x_n) = \sum_{i=1}^n\sum_{y\in S}
p(x_i,y)\left(a+b\sum_{j=1}^n I(y=x_j)\right) (f(x^{x_i,y})-f(x))
\ee
so
\[
 \loc^{a,b}_n = \loc^a_{1,n} + \loc^b_{2,n}
\]
where
\be\label{walkparta}
\loc^a_{1, n} f(x_1,\ldots,x_n) = a\sum_{i=1}^n\sum_{y\in S}
p(x_i,y)(f(x^{x_i,y})-f(x))
\ee
is the independent random walk part (random walks jumping at rate $a$)
and
\be\label{inparta}
\loc^b_{2,n} f(x_1,\ldots,x_n) = b\sum_{i=1}^n\sum_{j=1}^n
p(x_i,x_j)(f(x^{x_i,x_j})-f(x))
\ee
is the ``clumping'' part, i.e., when $b<0$ clumping is discouraged,
and $b>0$ clumping is favored.

We call $T^{a,b}_n(t)$ the semigroup on functions $f:S^n\to \R$
associated to the generator \eqref{abgen},
and $U^{a}_n(t)$ the semigroup of a system of
independent continuous-time random walkers (jumping at rate $a$), i.e.,
the semigroup associated to the generator $\loc^{a}_{1, n}$ in \eqref{walkparta}.
Notice that when $b<0$, $T_n^{a,b}(t)$ is not
always a Markov semigroup. However, for the applications of negative $b$,
we have in mind generalized exclusion process (see Section 6)
in which case $a/b$ is an integer and in this case $T^{a,b}_n(t)$ is
a Markov semigroup.

\bt\label{siplig}
Let $f:S^n\to\R$ be positive definite and symmetric. Then we have
for $b>0$
\be\label{ligin}
U^a_n (t) f\leq T^{a,b}_n (t) f
\ee
and
for $b<0$, if $(T^{a,b}(t))_{t\geq 0}$ is a Markov semigroup, we have
\be
U^a_n (t) f\geq T^{a,b}_n (t) f
\ee
\et
\bpr
The proof follows the proof in \cite{ligg}, but for the sake of self-constistency we prefer to give
it explicitely. Suppose $b>0$.

Start with the decomposition \eqref{finsipgen} and
use the symmetry of $p(x,y)$ and $f$ to write
\beq
&&(\loc^{a,b}_n f - \loc^{a}_{1,n} f)(x) = (\loc^{b}_{2,n} f)(x)
\nonumber\\
&= &
b\sum_{i=1}^n\sum_{j=1}^n
p(x_i,x_j)(f(x^{x_i,x_j})-f(x))
\nonumber\\
&=&
\frac{b}{2}
\sum_{i=1}^n\sum_{j=1}^n
p(x_i,x_j)(f(x^{x_i,x_j}) + f(x^{x_j,x_i})-2f(x))
\nonumber\\
&=&
\frac{b}{2}\sum_{i=1}^n\sum_{j=1}^n
p(x_i,x_j)
\nonumber\\
&\times &
\sum_{u,v}f(x_1,\ldots, x_{i-1},u, x_{i+1},\ldots, x_{j-1}, v, x_{j+1}, \ldots, x_n)
(\delta_{x_i,u}-\delta_{x_j,u})(\delta_{x_i,v}-\delta_{x_j,v})
\nonumber\\
&\geq & 0 \eeq
where in the last step we used that $f$ is positive definite.

Since
$U^a_n(t)$ is the semigroup of independent walks, it maps positive
definite functions into positive definite functions, and so we have
\[
 \left(\loc_n U^a_n(t) f - \loc^a_{1,n} U^a_n (t) f\right)= \loc^b_{2,n} U^a_n(t) f\geq 0
\]
We can then use the variation of constants formula
\be\label{theend}
 T^{a,b}_n (t)f - U^a_n(t)f = \int_0^t ds\ T^{a,b}_n (t-s) \left(\loc^b_{2,n} U^a_n(s) f\right)\geq 0
\ee
and remember that $T^{a,b}_n(t)$ is a Markov semigroup which therefore maps
non-negative functions into non-negative functions.

The proof for $b<0$, under the assumption that $T^{a,b}_n(t)$ is a Markov semigroup is identical.
\epr

\section{Stationary measures and self-duality for the $SIP(m)$}
The stationary measures of $\sip$ are product measures
of ``discrete gamma distributions''
\[
 \nu_\lambda (d\eta) = \otimes_{x\in S} \nu^m_\lambda (d\eta_x)
\]
where
for $n\geq 0$
\be\label{revprod} \nu^m_\lambda (n) = \frac{1}{Z_{\lambda
,m}}\frac{\lambda^n}{n!} \frac{\Gamma (\frac{m}2 +n)}{\Gamma
(\frac{m}2 )}, \qquad n\in\N \ee
with $0\leq \lambda < 1$ a parameter, $\Gamma(r)$  the gamma-function
and
\[
 Z_{\lambda,m} = \left(\frac{1}{1-\lambda}\right)^{m/2}
\]

Notice that for $m=2$, $\nu_\lambda^m$ is a geometric distribution
(starting from zero), i.e., $\nu_\lambda^2 (n) = \lambda^n
(1-\lambda), n\in\N$ and for $m/2$ an integer $\nu_\lambda^m$
is negative binomial distribution $NB(m/2,\lambda)$.
Moreover, the measures $\nu^m$ have the
following convolution property \be\label{koet}
\nu^m_\lambda*\nu^l_\lambda= \nu^{m+l}_\lambda \ee where $*$ denotes
convolution, i.e., a sample from $\nu^m_\lambda*\nu^l_\lambda$ is
obtained by site-wise addition of a sample from $\nu^m_\lambda$ and
an independent sample from $\nu^l_\lambda$.

The $SIP(m)$ process is self-dual \cite{gkrv} with duality functions
given by
$D(\xi,\eta) = \prod_x d(\xi_x,\eta_x)$, with
\be\label{vleermuis}
 d(k,l) = \frac{l!}{(l-k)!} \frac{\Gamma \left(\frac{m}2\right)}{\Gamma\left(\frac{m}2 +k\right)}
\ee
where $k \le l$. Self-duality means that
\be\label{sipself}
\E^{\sip}_\eta D (\xi,\eta_t) = \E^{\sip}_\xi D(\xi_t,\eta)
\ee
where $\E^{\sip}_\eta$ denotes expectation in the SIP process started from
the configuration $\eta$.

The relation between the polynomials $D$ and the measure
$\nu^m_\lambda$ reads \be\label{mpolrel} \int
D(\xi,\eta)\nu^m_\lambda (d\eta) =
\left(\frac{\lambda}{1-\lambda}\right)^{|\xi|} \ee as follows from a
simple computation using the definition of the $\Gamma$-function,
$\Gamma(r) =\int_0^\infty x^{r-1} e^{-x} dx$.

From conservation of particles in the dual process, we see
that self-duality and the relation \eqref{mpolrel} gives
stationarity of the measure $\nu_\la$.

The relation \eqref{mpolrel} can be generalized to ``local stationary measure'',
i.e. the product measures that are obtained from the
stationary measure \eqref{revprod} by allowing a site-dependent
parameter. More precisely, given
\[
 \laa : S\to [0,1)
\]
we define the local stationary measure associated
to the profile $\laa$ by
\be\label{profilemeas}
 \nu_{\laa} = \otimes_{x\in S} \nu^m_{\laa (x)} (d\eta_x)
\ee
For $x_1,\ldots,x_n\in S$ we denote
by $\sum_{i=1}^n \delta_{x_i}$ the particle
configuration $\xi\in \N^S$ obtained by putting a particles at
locations $x_i$, i.e., $\xi (x)= \sum_{i=1}^n I(x_i=x)$.
We then have the following relation between the duality functions
and the local stationary measures
\be\label{estrel}
\int D \left(\six,\eta\right) \nu_{\laa} (d\eta) = \prod_{i=1}^n \rho (x_i)
\ee
where
\[
 \rho (x_i) = \frac{\lambda (x_i)}{1-\lambda (x_i)}
\]
For a constant profile $\laa(x) =\lambda\;,\forall x\in S$, we recover \eqref{mpolrel}.

By Lemma \ref{vanishlem}
below, in the case $S=\Zd$ and translation
invariant finite-range $p(x,y)$, any number
of dual particles in the $\sip$ will eventually diffuse
away to infinity.
From that it is easy to deduce that
the measures $\nu_\lambda$ are extremal invariant.
To see
this, we denote for two finite particle configurations
$\xi\perp\xi'$, if their supports are disjoint, i.e., there are no
site $x\in S$ where there are $\xi$ and $\xi'$ particles. If
$\xi\perp\xi'$ then $D(\xi+\xi',\eta) = D(\xi,\eta) D(\xi',\eta)$.
Since at large $t>0$, assumption (A2) implies that, in the SIP
started with a finite number of particles, particles are with
probability close to one at different locations (see Lemma
\ref{vanishlem} for a proof of this), we have that for $\xi'$ a
fixed configuration, the event $\xi_t\perp\xi'$ has probability
close to one as $t\to\infty$. Therefore
\beq &&\lim_{t\to\infty}\int
\E^{SIP(m)}_\eta\left(D(\xi,\eta_t)\right) D(\xi',\eta) \nu_\lambda
(d\eta)
\nonumber\\
&=& \lim_{t\to\infty}\E^{SIP(m)}_\xi \int D(\xi_t, \eta)
D(\xi',\eta) \nu_\lambda (d\eta)
\nonumber\\
&=& \lim_{t\to\infty}\E^{SIP(m)}_\xi \int D(\xi_t, \eta) D(\xi',\eta) I
(\xi_t \perp \xi')  \nu_\lambda (d\eta)
\nonumber\\
&=&
\lim_{t\to\infty}\rho_\lambda^{|\xi_t| + |\xi'_t|}
\nonumber\\
&=&
\rho_\lambda^{|\xi| + |\xi'|}
\nonumber\\
&=& \int D(\xi,\eta)\nu_\lambda (d\eta)\int D(\xi',\eta) \nu_\lambda
(d\eta) \eeq
which shows that time-dependent correlations of (linear combinations of) $D(\xi,\cdot)$
polynomials decay in the course of time to zero, and hence, by standard arguments,
$\nu_\lambda$ is mixing and thus ergodic.

\section{Correlation inequalities in the $\sip$}\label{corre}
For a probability measure $\mu$ on the configuration space
$\N^S$, we denote its ``duality moment function''
$K_{\mu}:S^n\to\R$ by
\be\label{fourier}
 K_{{\mu}} (x_1,\ldots,x_n)
= \int D \left(\sum_{i=1}^n \delta_{x_i} ,\eta\right)\mu (d\eta)
\ee
If $\mu=\nu_{\laa}$ is a local stationary measure with profile $\laa$, then
\be\label{locfou}
K_{\nu_{\laa}} (x_1,\ldots, x_n) = \prod_{i=1}^n \rho (x_i)
\ee
which is clearly positive definite and symmetric.
We can therefore apply Theorem \ref{siplig}
and obtain the following result.
\bp\label{sipcorprop} For all $t\geq 0$, for all
profiles $\laa: S\to [0,1)$ and for all
$x_1,\ldots,x_n\in S$ we have
\be\label{compsip}
K_{\nu_{\laa} S_t} (x_1,\ldots,x_n) \geq \prod_{i=1}^n K_{\nu_{\laa} S_t} (x_i)
\ee
where $S_t$ denotes the semigroup of
the $\sip$ process.
In particular, when the $\sip$ is started from $\nu_{\laa}$,
the random variables $ \{ \eta_t (x), x\in S\}$ are positively
correlated, i.e., for $(x,y)\in S\times S$
$$
\int
\E^{\sip}_{\eta}\left(\eta_t(x)\eta_t(y)\right)\nu_{\laa}(d\eta)
\ge \int\E^{\sip}_{\eta}\left(\eta_t(x)\right)\nu_{\laa}(d\eta)
\int\E^{\sip}_{\eta}\left(\eta_t(y)\right)\nu_{\laa}(d\eta)
$$
\ep
\bpr Denote by $\E^{\sip}_{x_1,\ldots,x_n}$ expectation in the $\sip$
process started with $n$ particles at positions $(x_1,\ldots,x_n)$,
by $\E^{IRW(m)}$ expectation in the process of independent random
walkers (jumping at rate $2m$) and $\E^{RW(m)}$ a single random walker expectation.
We then have the following chain of inequalities,
which is obtained by using sequentially the following:
self-duality property
(\ref{sipself}), the comparison inequality (\ref{ligin}), the
relation between the measure $\nu_{\laa}$ and the duality
function $D$ (\ref{estrel}), the independence between random
walkers, the fact that a single SIP particle moves as a continuous
time random walk, and finally again self-duality
(\ref{sipself})
\beq
&&\int \E^{\sip}_\eta D\left(\sum_{i=1}^n
\delta_{x_i},\eta_t\right) \nu_{\overline{\lambda}} (d\eta)
\nonumber\\
&= &
\E^{\sip}_{x_1,\ldots,x_n}\int D\left(\sum_{i=1}^n \delta_{X_i(t)},\eta\right) \nu_{\laa} (d\eta)
\nonumber\\
&\geq &
\E^{IRW(m)}_{x_1,\ldots,x_n}\int D\left(\sum_{i=1}^n \delta_{X_i(t)},\eta\right) \nu_{\laa} (d\eta)
\nonumber\\
&=&
\E^{IRW(m)}_{x_1,\ldots,x_n} \left(\prod_{i=1}^n \rho (X_i(t))\right)
\nonumber\\
&=&
\prod_{i=1}^n \E^{RW(m)}_{x_i} \rho (X_i(t))
\nonumber\\
&=&
\prod_{i=1}^n \int \E^{\sip}_{x_i} \left(D(\delta_{X_i(t)}, \eta)\right)\nu_{\laa} (d\eta)
\nonumber\\
&=& \prod_{i=1}^n \int\E^{\sip}_{\eta}
\left(D(\delta_{x_i},\eta_t)\right)\nu_{\laa} (d\eta)
\eeq
\epr

This proposition shows that
starting from a local stationary measure
$\nu_{\laa}$, the density profile
$\rho_t (x)= \E^{RW(m)}_x \rho_t (x)$ predicts (by duality)
correctly the density at time $t>0$ but the
true measure at time $t>0$, $\nu_{\laa} S_t$, lies above
(in the sense of expectations of $D$-functions)
the product measure with density profile $\rho_t (x)$.

From the analogy with the SEP emphasized above, one could think that
\eqref{compsip} extends to the case when the SIP
process is started from a general product measure. However,
for general probability measures $\mu$ on $\Omega$, the duality moment function
$K_{\mu}: S^n \to R$ defined in \eqref{fourier}
is not
necessarily positive definite (as is the case for the special
product measures $\nu_{\laa}$), since we do not have the equality
$ D\left(\sum_{i=1}^n \delta_{x_i},\eta\right) =
\prod_{i=1}^n  D(\delta_{x_i},\eta)$ in general. Notice that this
problem does not appear in the context of the standard SEP, as for that
model, the self-duality functions are
\[
 D_{SEP}\left(\sum_{i=1}^n \delta_{x_i},\eta\right) = \prod_{i=1}^n \eta_{x_i}= \prod_{i=1}^n
D_{SEP}\left( \delta_{x_i},\eta\right)
\]
and hence automatically, for any measure $\mu$, the function $K_{\mu}$
is positive definite in that model.

If however {\em all $x_i$ are different}, we have $D(\sum_{i=1}^n
\delta_{x_i},\eta)= \prod_{i=1}^n D( \delta_{x_i},\eta)$. For
every probability measure $\mu$ on $\Omega$, the function $\Psi_\mu:
S^n\to\R$ defined by \be\label{croki1}
 \Psi_\mu (x_1,\ldots,x_n) = \int \prod_{i=1}^n D(\delta_{x_i},\eta)
 \mu(d\eta)
\ee is clearly positive definite. This, together with the fact that
under assumption (A2), a finite number of $\sip$ particles diffuse and
therefore eventually will be typically at different positions,
suggests that in a stationary measure, the variables $\eta_{x_i}$
are positively correlated.

To state this result we introduce the class of probability measures
with uniform finite moments \be\label{classpr} \caP_f =: \{ \mu:
\forall n\in \N,\ \sup_{|\xi|=n} \int D(\xi,\eta) \mu (d\eta)=:
M^n_\mu <\infty\} \ee For a sequence of measures $\mu_n\in \caP_f$,
and $\mu\in \caP_f$, we define that $\mu_n\to \mu$ if for all $\xi$
finite particle configuration,
\[
 \lim_{n\to\infty} \int D(\xi,\eta) \mu_n (d\eta) =\int D(\xi,\eta) \mu (d\eta)
\]

We can then formulate our next result. \bp\label{geninvprop} Assume
(A1) and (A2). Let $\nu\in\caP_f$ be a product
measure. Let $S(t)$ denote the semigroup of the \sip. Suppose that
\be\label{subseqdef}
\mu= \lim_{n\to\infty} \nu S(t_n)
\ee for a subsequence
$t_n\uparrow\infty$. Then we have $\mu\in \caP_f$, $\mu$ is
invariant and \be\label{corsipin} K_{\mu} (x_1,\ldots,x_n)\geq
\prod_{i=1}^n K_{\mu} (x_i) \ee \ep \bpr First, by duality we
have, referring to the definition of $\caP_f$, for all $t>0$,
\[
 \int \E^{\sip}_\eta D(\xi,\eta_t ) \nu (d\eta) = \E^{\sip}_{\xi}\int D(\xi_t,\eta) \nu (d\eta)
\leq M^{|\xi|}_\nu <\infty
\]
which shows that both $\nu S(t_n)$ and $ \mu$ are elements of
$\caP_f$. The invariance of $\mu$ follows from duality,
$\nu\in\caP_f$ and Lemma 1.26 in \cite{ligg}, chapter V.

To proceed with the proof of the proposition, we start with the
following lemma, which ensures that, under condition (A2), any
number of $\sip$ particles will eventually be at different locations.
\bl\label{vanishlem} Assume (A2). Start the finite $\sip$ with
particles at locations $\{ x_1,\ldots, x_n\}$, then
\be\label{sipdif} \lim_{t\to\infty}\pee^{\sip}_{x_1,\ldots,x_n}
\left(\exists i\not=j: X_i(t)=X_j(t)\right)=0 \ee
\el

\bpr
We give the proof for $m=1$. The general
case is a straightforward extension.
Put
$\eta:= \sum_{i=1}^n \delta_{x_i}$. Using self-duality we can write
\beq\label{lalala} \pee^{SIP(1)}_\eta \left(\exists i\not= j:
X_i(t)=X_j(t)\right) &\leq & \sum_z \pee^{SIP(1)}_\eta \left(\eta^2_t
(z)-\eta_t (z) > 1\right)
\nonumber\\
&\leq &
\sum_z \E^{SIP(1)}_\eta (\eta^2_t (z)-\eta_t (z))
\nonumber\\
&=&
\frac34\sum_z \E^{SIP(1)}_\eta \left(D(2\delta_z, \eta_t)\right)
\nonumber\\
&=&
\frac34 \sum_{z} \E^{SIP(1)}_{z,z}\left( D(\delta_{X_t}+ \delta_{Y_t}, \eta)\right)
\nonumber\\
&\leq &
3 \sum_{z} \E^{SIP(1)}_{z,z} (\eta (X_t)\eta (Y_t))
\nonumber\\
&= &
3 \sum_{z}\sum_{i,j=1}^n \E^{SIP(1)}_{z,z} \left(I(X_t= x_i) I(Y_t= x_j)\right)
\nonumber\\
&\leq & 3n^2 \sup_{x,y} \pee^{\sip}_{x,y} (X_t=Y_t) \eeq where in the
last step we used the symmetry of the transition probabilities of
the $SIP(1)$ (with two particles). \epr We now proceed with the proof of
the proposition. For $x_1,\ldots,x_n \in S$ we define
\be\label{loop} \left|D(\sum_{i=1}^n \delta_{x_i},\eta)-
\prod_{i=1}^n D(\delta_{x_i},\eta)\right|=
\Delta(x_1,\ldots,x_n,\eta) \ee We have that $\Delta
(x_1,\ldots,x_n,\eta)=0$ if all $x_i$ are different, i.e., if $|\{
x_1,\ldots,x_n\}|=n$.
Since by assumption (A2) and
Lemma \ref{vanishlem}, the probability that two $\sip$ walkers out of a
finite number $n$ of them occupy the same position, i.e.
$X_i(t)=X_j(t)$ for some $i \not= j$, vanishes in the limit
$t\to\infty$, we conclude, using $\nu\in\ \caP_f$, for any
$x_1,\ldots, x_n\in S$, \be\label{bambam}
 \lim_{t\to\infty}\int \E^{\sip}_{x_1,\ldots,x_n}\Delta(X_1(t),\ldots,X_n(t),\eta)\nu(d\eta)=0
\ee
Moreover from the comparison inequality
(\ref{ligin}) we have, using the notation \eqref{croki1}
\beq\label{gier} \E^{\sip}_{x_1,\ldots,x_n} \Psi_\nu
(X_1(t),\ldots,X_n (t)) &\geq & \E^{IRW(m)}_{x_1,\ldots,x_n}\Psi_\nu
(X_1(t),\ldots,X_n (t))
\nonumber\\
&=&
\E^{IRW(m)}_{x_1,\ldots,x_n}\int \prod_{i=1}^n D\left( \delta_{X_i (t)},\eta\right) \nu (d\eta)
\nonumber\\
&=& \prod_{i=1}^n\E_{x_i}^{RW(m)}\int D(\delta_{X_i (t)},\eta) \nu
(d\eta)+ \epsi (t) \eeq where $\epsi (t)\to 0$ as $t\to\infty$ by
assumption (A1), i.e., for large $t>0$, independent random walkers
are at different locations with probability close to one. Therefore,
using the definition (\ref{subseqdef}), the self-duality property
(\ref{sipself}), the equation \eqref{bambam}, the equation
\eqref{gier}, and taking limits along the subsequence $t_n$ we have
\beq K_{\mu} (x_1,\ldots, x_n) &=&\lim_{t\to\infty} \int
\E^{\sip}_{\eta} D\left(\sum_{i=1}^n \delta_{x_i},\eta_t\right)\nu
(d\eta)
\nonumber\\
&=& \lim_{t\to\infty} \int \E^{\sip}_{x_1,\ldots,x_n}
D\left(\sum_{i=1}^n \delta_{X_i (t)},\eta\right)\nu (d\eta)
\nonumber\\
&=& \lim_{t\to\infty} \E^{\sip}_{x_1,\ldots,x_n}\Psi_\nu (X_1 (t),\ldots, X_n (t))
\nonumber\\
&\geq &
\lim_{t\to\infty} \prod_{i=1}^n\E_{x_i}^{RW(m)}\int D(\delta_{X_i (t)},\eta) \nu (d\eta)
\nonumber\\
&=&
\prod_{i=1}^n K_{\mu } (x_i)
\eeq
\epr

\section{Correlation inequalities in the SEP$(n)$}
We now consider the application of the generalized Liggett inequality
for negative $b$.
The $SEP(n)$ is the Markov process on $\Omega = \{ 0,1,\ldots, n\}^S$ with
generator
\be\label{sepk}
L f(\eta) = \sum_{x,y\in S} \eta(x) (n-\eta(y)) p(x,y) \left(f(\eta^{xy}) -f(\eta)\right)
\ee
The stationary measures of this process are products of binomial distributions, i.e.,
for $\rho\in [0,1]$,
\be\label{binomial}
\nu_\rho = \otimes_{x\in S} Bin( n,\rho)
\ee
Similar to the case of the inclusion process, for a
profile
$\overline{\rho} : S\to [0,1]$
we define the local stationary measure
\[
 \nu_{\overline{\rho} }= \otimes_{x\in S} Bin (n,\overline{\rho}(x))
\]
The duality functions for self-duality are given by (see \cite{gkrv})
\be\label{sepdual}
D(\xi,\eta) = \prod_x d(\xi_x,\eta_x)
\ee
for $\xi\in\Omega$ a configuration with finitely many particles
(at most $n$ per site) and
with
\be\label{binpol}
d(k,l) = \frac{{l\choose k} }{{n\choose k}}
\ee
The relation between the duality functions and the local
stationary measures is, as usual,
i.e., for $\xi = \sum_{i=1}^n \delta_{x_i}\in \Omega$ (i.e., at most
$n$ particles per site), and $\roo$ a profile:
\be\label{dualrelsep}
\int D(\xi, \eta) \nu_{\roo} (d\eta) = \prod_{i=1}^n \rho (x_i)
\ee
We define, for a probability measure $\mu$ on $\Omega$, its duality moment function
\be\label{sepfourier}
K_{\mu} (x_1,\ldots,x_n) = \int D\left(\sum_{i=1}^n \delta_{x_i},\eta\right) \mu(d\eta)
\ee

The following proposition is then the analogue of
Proposition \ref{sipcorprop} in this context (with inequality in
the other direction since $b<0$).
\bp\label{sepcorprop}
For $\roo: S\to [0,1]$ a density profile and
$ t>0$,
\be\label{hatsep}
K_{\nu_{\roo} S_t} (x_1,\ldots,x_n) \leq \prod_{i=1}^n K_{\nu_{\roo} S_t} (x_i)
\ee
In particular, for starting from
$\nu_{\roo}$, the variables
$\{ \eta_t(x): x\in S\}$ are negatively correlated.
\ep
\section{Correlation inequalities for some interacting diffusions}

\subsection{The Brownian Momentum Process}
The Brownian momentum process is a system of interacting diffusions,
initially introduced as a model of heat conduction in \cite{gk}, and
analyzed via duality in \cite{gkr}. It is defined as a Markov
process on $X=\R^S$ via the formal generator on local functions:
\be\label{bmpgen} L_{BMP} f(\eta) = \left(\sum_{x,y\in S} p(x,y)
\left(\eta_x\frac{\partial}{\partial \eta_{y}}
-\eta_x\frac{\partial}{\partial \eta_{y}}\right)^2\right) f (\eta)
\ee The variable $\eta_x$ has to be thought of as momentum of an
``oscillator'' associated to the site $x\in S$. The local kinetic energy
$\eta_x^2$ has to be thought of as the analogue of the number of
particles at site $x$ in the $\sip$ with $m=1$. The expectation of $\eta_x^2$ is interpreted
as the local temperature at $x$.

Defining the polynomials
\[
D(n,z)= \frac{z^{2n}}{(2n-1)!!}
\]
we have the duality function $D(\xi,\cdot)$ defined
on $X$ and indexed by finite particle configurations
$\xi\in \N^S, \sum_x \xi_x <\infty$:
\be
D(\xi,\eta) = \prod_{x\in S} D(\xi_x,\eta_x)
\ee
In \cite{gkr}, \cite{gkrv}, we proved
the duality relation
\be\label{duality}
\E^{BMP}_\eta \left(D(\xi,\eta_t)\right)= \E^{\siP}_\xi \left(D(\xi_t,\eta)\right)
\ee
As before, for $x_1,\ldots,x_n\in S$ we denote by $\sum_{i=1}^n \delta_{x_i}$ the particle
configuration obtained by putting a particle at each $x_i$.

Let $\mu$ be a product of Gaussian measures on $X$, with
site-dependent variance, i.e., for a function $\rho: S\to
[0,\infty)$, we define \be\label{gauss} \mu_\rho = \otimes_{x\in S}
\nu_{\rho(x)} (d\eta_x) \ee where
\[
\nu_{\rho(x)}(d\eta_x)  = \frac{e^{-\eta_x^2/2\rho(x)}}{\sqrt{2\pi
\rho(x)}} d\eta_x
\]
is the Gaussian measure on $\R$ with mean zero and variance $\rho
(x)$. Then we have \be \int D\left(\sum_{i=1}^n \delta_{x_i},\eta\right)
\mu_\rho (d\eta) = \prod_{i=1}^n \rho (x_i) \ee From this
expression, it is obvious that the map \be\label{psimap} S^n \to \R
: (x_1,\ldots,x_n)\mapsto \int D\left(\sum_{i=1}^n \delta_{x_i},\eta\right)
\mu_\rho (d\eta) \ee is positive definite. Therefore, combining
the duality property between BMP process and $\sip$
process, \eqref{duality}, with Theorem \ref{siplig} we have the
inequality
\beq &&\int \E^{BMP}_\eta D\left(\sum_{i=1}^n
\delta_{x_i},\eta_t\right) \mu_\rho (d\eta)
\nonumber\\
&= & \E^{\siP}_{x_1,\ldots,x_n}\int D\left(\sum_{i=1}^n
\delta_{X_i(t)},\eta\right) \mu_\rho (d\eta)
\nonumber\\
&\geq & \E^{IRW(m)}_{x_1,\ldots,x_n}\int D\left(\sum_{i=1}^n
\delta_{X_i(t)},\eta\right) \mu_\rho (d\eta)
\nonumber\\
&=& \E^{IRW(m)}_{x_1,\ldots,x_n} \left(\prod_{i=1}^n  \int
D\left(\delta_{X_i(t)},\eta\right) \mu_\rho (d\eta) \right)
\nonumber\\
&=&
\E^{IRW(m)}_{x_1,\ldots,x_n} \left(\prod_{i=1}^n \rho (X_i(t))\right)
\nonumber\\
&=& \prod_{i=1}^n \E^{RW(m)}_{x_i} \rho (X_i(t)) \nonumber\\
& = & \prod_{i=1}^n \int
\E^{\siP}_{x_i}\left(D\left(\delta_{X_i(t)},\eta\right)\right) \mu_\rho (d\eta) \nonumber\\
& = & \prod_{i=1}^n \int
\E^{BMP}_{\eta}\left(D\left(\delta_{x_i},\eta_t\right)\right) \mu_\rho (d\eta) \eeq which
is the analogue of Proposition \ref{sipcorprop} for the BMP process.

In words, it means that the ``non-equilibrium temperature profile''
is above the temperature profile predicted from the discrete diffusion
equation.
It also implies that
the variables $\{\eta_x^2: x\in S\}$ are positively correlated under
the measure $(\mu_\rho)_t$ for all choices of $\rho$, $t>0$.

More precisely, if we denote
\[
 \rho_t (x)= \E^{RW(m)}_{x} \rho (X_t)
\]
then we have that $\eta_x^2$ at time $t$ has
expectation $\rho_t(x)$ when the starting measure
is $\mu_\rho$ (since a single particle in the $\siP$ moves as
a continuous time random walk). The correlation inequality
for the BMP which we just derived shows that the true
measure at time $t>0$ when started from a product of Gaussian measures
lies stochastically above the
Gaussian product measure with mean zero and variance
$\rho_t (x)$.

Similarly, we obtain an analogous correlation inequality
for the BMP
for a measure obtained as a limit of
product measures.
We define
\[
 \caP_f (X) = \{ \mu: \forall n\in \N: \sup_{|\xi|=n} \int D(\xi,\eta) \mu (d\eta) <\infty\}
\]

\bp Assume (A1) and (A2). Suppose $\nu\in\caP_f
(X)$ is a product measure and $\mu$ is a limit
point of the set $\{ \nu S(t): t\geq 0\}$, where $S(t)$ denotes the
semigroup of the BMP process. Then we have the inequality
\[
 K_{\mu} (x_1,\ldots,x_n) \geq \prod_{i=1}^n K_{\mu} (x_i)
\]
\ep

\subsection{The Brownian Energy Process}
The Brownian energy process with parameter
$m>0$ (notation $BEP(m)$) is introduced in
\cite{gkrv} as the process on state space $X= [0,\infty)^S$,
with generator
\be\label{bepm}
L  = \sum_{x,y\in S} p(x,y) L^m_{xy}
\ee
with
\be\label{bepmxy}
L^m_{xy} f(\eta)  = 4 \eta_x\eta_y
\left(\frac{\partial}{\partial\eta_x}-
\frac{\partial}{\partial\eta_x}\right)^2 f(\eta) - 2m (\eta_x-\eta_y)
\left(\frac{\partial}{\partial\eta_x}-
\frac{\partial}{\partial\eta_x}\right) f(\eta)
\ee
This process is dual to the $\sip$ in the following sense.
Define, for $\xi\in \N^S$ a finite particle configuration, and
$\eta\in X$ the polynomials
\be\label{beppol}
D(\xi,\eta) = \prod_{x\in S} d(\xi_x,\eta_x)
\ee
with, for $k\in\N, y\in [0,\infty)$
\be\label{bepsingle}
d(k,y) = y^{k} \frac{\Gamma \left(\frac{m}{2}\right)}{2^k \Gamma \left(\frac{m}{2}+k\right)}
\ee
then we have
\be\label{bepsipdual}
\E^{BEP(m)}_\eta D(\xi,\eta_t) = \E^{\sip}_\xi D(\xi_t,\eta)
\ee
As a consequence, extremal invariant measure of the $BEP(m)$ are products
of $\Gamma$-distributions with shape parameters $m/2$ and scale parameter $\theta>0$:
\be\label{bepstat}
\nu_{\theta} (d\eta) = \otimes_{x\in S} \nu_\theta (d\eta_x)
\ee
with
\be\label{bepgamma}
\nu_\theta (dz) = \frac{1}{\theta^{m/2} \Gamma \left(\frac{m}{2}\right)}z^{\frac{m}{2} -1} e^{-z/\theta}
\ee
Similarly we define the local stationary measures
\be\label{locbepstat}
\nu_{\overline{\theta}} =
\otimes_{x\in S} \nu_{\overline{\theta} (x)} (d\eta_x)
\ee
with $\overline{\theta}: S\to [0,\infty)$, and the duality moment function
of a probability measure $\mu$ on $X$:
\be\label{dualmomentbep}
K_\mu (x_1,\ldots,x_n) =  \int D\left(\sum_{i=1}^n \delta_{x_i}, \eta\right) \mu (d\eta)
\ee
As a consequence of the correlation inequalities derived for
the $\sip$, we derive the following.
\bp\label{bepprop}
\ben
\item
For all $\overline{\theta}: S\to [0,\infty)$, $t>0$, and $x_1,\ldots, x_n \in S$ we have
\be\label{bepcor}
K_{\nu_{\overline{\theta}}S^{BEP(m)}_t} (x_1,\ldots,x_n)\geq \prod_{i=1}^n K_{\nu_{\overline{\theta}}S^{BEP(m)}_t}\mu (x_i)
\ee
\item If for some product measure $\nu$ on $X$ with finite moments,
and a sequence of $t_n\uparrow\infty$
the limit
\[
\mu = \lim_{n\to\infty} \nu S^{BEP(m)}(t_n)
\]
exists, then
\be\label{statbepcor}
K_\mu (x_1,\ldots,x_n) \geq \prod_{i=1}^n K_\mu (x_i)
\ee
\een
\ep

\section{The boundary driven $\sip$}

In this section we consider the non-equilibrium
one-dimensional model that is obtained by considering particle
reservoirs attached to the first and last sites of the chain. We
will show that, if one requires reversibility w.r.t. the measure
$\nu_{\lambda}^m$ {\em and} duality with absorbing boundaries, this
uniquely fixes the birth and death rates at the boundaries.
\subsection{Duality for the the boundary driven $\sip$}
The generator of the boundary driven $\sip$ on a chain
$\{1,\ldots,N\}$ driven at the end points, reads \be\label{bdsip} \loc
= \loc_1 +\loc_N + \loc_{bulk} \ee where $\loc_{bulk}$ denotes the
$\sip$ generator, with nearest neighbor random walk as underlying
kernel, i.e., \be\label{bulkgen}
 \loc_{bulk} f(\eta) = \sum_{x\in \{1,\ldots,N-1\}}  2\eta_x (m+2\eta_{x+1}) \left(f(\eta^{x,x+1}) -f(\eta)\right)
+ 2\eta_{x+1} (m+2\eta_x) \left(f(\eta^{x+1,x}) -f(\eta)\right)
\ee
and
where $\loc_1,\loc_N$ are birth and death processes
on the first and $N$-th variable respectively, i.e.,
\[
 \loc_{1} f(\eta) = d_L(\eta_1)(f(\eta-\delta_1)- f(\eta))  + b_L(\eta_1) (f(\eta+\delta_1)-f(\eta))
\]
and
\[
 \loc_{N} f(\eta) = d_R(\eta_N)(f(\eta-\delta_N)- f(\eta))  + b_R(\eta_N) (f(\eta+\delta_N)-f(\eta))
\]
These generators model contact with respectively the left and right particle reservoir.

The rates $d_L,b_L,d_R,b_R$ are chosen such that detailed balance is
satisfied w.r.t.\ the measure $\nu^m_\lambda$, with
$\lambda=\lambda_L$ for $d_L, b_L$, and $\lambda=\lambda_R$ for
$d_R, b_R$. More precisely, this means that these rates satisfy
\be\label{ciccio}
 b_\alpha (k) \nu^m_{\lambda_{\alpha}} (k) = d_{\alpha} (k+1) \nu^m_{\lambda_{\alpha}} (k+1)
\ee for $\alpha\in \{ L,R\}$.

To state our duality result, we consider functions $\caD (\xi,\eta)$
indexed by particle configurations $\xi$ on $ \{0,\ldots, N+1\}$
defined by \be\label{bddual} \caD (\xi,\eta)= \rho_L^{|\xi_0|}
D(\xi_{\{1,\ldots,N\}},\eta) \rho_R^{|\xi_{N+1}|} \ee where
$\rho_{\alpha} = \rho_{\lambda_{\alpha}} =
\lambda_{\alpha}/(1-\lambda_{\alpha})$, and where we remember that
\[
 D(k,n) = \frac{n!}{(n-k)!} \frac{ \Gamma \left(\frac{m}2\right)}{\Gamma \left(\frac{m}2+k\right)}
\]
is the duality function for the $\sip$.
I.e., for the ``normal'' sites $\{1,\ldots,N\}$ we simply
have the old duality functions, and for the ``extra added'' sites
$\{0, N+1\}$ we have the expectation of the duality function
over the measure $\nu^m_\lambda$.

We now want duality to hold with duality functions
$\caD$, and with a dual process that behaves
in the bulk as the $\sip$, and which
has {\em absorbing boundaries} at $\{0, N+1\}$.
More precisely, we want the generator of the dual process
to be
\be\label{kauw}
\hat{\loc} =\loc_{bulk} + \hat{\loc}_1 + \hat{\loc_N}
\ee
with $\loc_{bulk}$ given by \eqref{bulkgen}, and
\[
 \hat{\loc}_1 f(\xi) = \xi_1 \left(f(\xi^{1,0})- f(\xi)\right)
\]
\[
 \hat{\loc}_N f(\xi) = \xi_N \left(f(\xi^{N,N+1})- f(\xi)\right)
\]
for $\xi \in \N^{\{0,1\ldots, N+1\}}$.
The duality relation then reads, as usual,
\be\label{merel}
 \left(\loc \caD (\xi,\cdot )\right) (\eta) = \left(\hat{\loc} \caD (\cdot,\eta)\right) (\xi)
\ee Since self-duality is satisfied for the bulk generator with the
choice \eqref{bddual}, i.e., since
\[
 \left(\loc_{bulk} \caD (\xi,\cdot )\right) (\eta) = \left({\loc}_{bulk} \caD (\cdot,\eta)\right) (\xi)
\]
\eqref{merel} will be satisfied if we have the following relations
at the boundaries: for all $k\leq n$:
\beq\label{moervos1} && b_{\alpha}(n) (D(k,n+1)- D(k,n)) +
d_{\alpha}(n) (D(k,n-1)-D(k,n))
\nonumber\\
&=& k(D(k-1,n)\rho_{\alpha} - D(k,n)) \eeq
where
$\alpha\in \{ L,R\}$.

From detailed balance (\ref{ciccio}) we obtain \be\label{raterel}
d_{\alpha}(n) = \frac1\lambda_{\alpha} \left( \frac{n}{\frac{m}2 +
n-1}\right) b_{\alpha}(n-1) \ee Working out \eqref{moervos1} gives,
using \eqref{vleermuis}, \beq\label{colibri} &&b_{\alpha}(n) \left(
\frac{n+1}{n+1-k} -1\right) + d_{\alpha}(n) \left( \frac{n-k}{n}
-1\right)
\nonumber\\
&=& k \left( \frac{ \left(\frac{m}2 +
k-1\right)\rho_{\alpha}}{n-k+1} -1\right) \eeq which simplifies to
\be\label{pinguin}
 \frac{b_{\alpha}(n)}{n+1-k} - \frac{d_{\alpha}(n)}{n}
= \left( \frac{ \left(\frac{m}2 + k-1\right)\rho_{\alpha}}{n-k+1}
-1\right) \ee Choosing \be\label{death} d_{\alpha}(n) =
\frac{n}{1-\lambda_{\alpha}} \ee and by the detailed balance
condition (\ref{raterel}), \be\label{birth} b_{\alpha}(n) = \left(
\frac{m}2 + n\right) \frac{\lambda_{\alpha}}{1-\lambda_{\alpha}} \ee
it is then an easy
computation to see that \eqref{moervos1} is satisfied with the
choices \eqref{death}, \eqref{birth}. Indeed, \eqref{pinguin}
reduces to the simple identity
\[
 \left(\frac{m}2 +n\right) \left(\frac{\lambda}{1-\lambda}\right)\frac{1}{n+1-k} -\frac1{1-\lambda}
= \frac{\frac{m}2 +k-1}{n+1-k}\left(\frac{\lambda}{1-\lambda}\right) -1
\]
We remark that the requirement of detailed balance alone is not
sufficient to fix the rates uniquely. However, the additional
duality constraint \eqref{moervos1} does fix the rates to the unique
expression given by (\ref{death}) and (\ref{birth}).

As a consequence of duality with duality functions
\eqref{bddual}, we have that the boundary driven
$\sip$ with generator \eqref{bdsip} has
a unique stationary measure $\mu_{L,R}$ for which
expectations of the polynomials $D(\xi,\eta)$ are given
in terms of absorption probabilities:
\beq\label{densprof}
\int D(\xi,\eta) \mu_{L,R} (d\eta )
&=& \lim_{t\to\infty} \E_\eta\caD (\xi,\eta_t)
\nonumber\\
&=&
\lim_{t\to\infty}  \hat{\E}_\xi\caD (\xi_t,\eta)
\nonumber\\
&=&
\sum_{k,l: k+l=|\xi|} \rho_L^k\rho_R^l \hat{\pee}_\xi \left(\xi_\infty = k \delta_0 + l\delta_{N+1}\right)
\eeq
Here, $\hat{\E}_\xi$ denotes expectation in the dual
process (which is absorbing at $\{0, N+1\}$) starting from
$\xi$.
In particular, since a single $\sip $ particle performs continuous
time simple random walk (at rate $2m$)
we have a linear density profile, i.e.,
\be\label{dens}
\int D(\delta_i,\eta) \mu_{L,R} (d\eta ) = \rho_L \left(1-\frac{i}{N+1}\right) + \rho_R \frac{i}{N+1}
\ee

\subsection{Correlation inequality for the boundary driven $\sip$}
For $x_1,\ldots,x_n \in \{ 1,\ldots, N\}$ let us denote by $(X_1
(t),\ldots, X_n (t))$ the positions of  particles at time $t$
evolving according to the $\sip$ with absorbing boundary sites at $ \{0
,N+1\}$, i.e., according to the generator \eqref{kauw}, and
initially at positions $x_1,\ldots,x_n$. Let $(Y_1 (t), \ldots, Y_n
(t))$ denote the positions at time $t$ of independent random walkers
(jumping at rate $2m$) absorbed (at rate 1) at $\{0, N+1\}$, initially
at positions $x_1,\ldots, x_n$. Since the absorption parts of the
generators of $(X_1 (t),\ldots, X_n (t))$ and $(Y_1 (t), \ldots, Y_n
(t))$ are the same, we have the same inequality for expectations of
positive definite functions as in Theorem \ref{siplig}. Therefore,
we have the following result on positivity of correlations in the
stationary state. This has once more to be compared to the analogous
situation of the boundary driven exclusion process, where the stationary
covariances of site-occupations are negative. \bp\label{kraai} Let
$\mu_{L,R}$ denote the unique stationary measure of the process with
generator \eqref{bdsip}. Let $x_1,\ldots,x_n\in \{ 1,\ldots, N\}$,
then we have \be\label{noneqcor} \int D\left( \sum_{i=1}^n
\delta_{x_i}, \eta\right) \mu_{L,R} (d\eta) \geq \prod_{i=1}^n\int
D(\delta_{x_i},\eta) \mu_{L,R} (d\eta) \ee In particular, $\eta_x,
x\in \{1,\ldots,N\}$ are positively correlated under the measure
$\mu_{L,R}$. \ep

\bpr Start from the measure $\nu_\lambda^m$. Define the map
$\{0,\ldots, N+1\}^n\to \R$: \be (x_1,\ldots,x_n) \mapsto \int
\caD\left( \sum_{i=1}^n \delta_{x_i},
\eta\right)\nu_\lambda^m (d\eta) = \prod_{i=1}^n
\rho(x_i) \ee where $\rho(x)= \frac{\lambda}{1-\lambda}$ for $x\in \{
1,\ldots, N\}$ and $\rho(0)= \rho_L, \rho (N+1) =\rho_R$. This is
clearly positive definite. Therefore, for $x_1,\ldots,x_n \in
\{1,\ldots,N\}$, we have \beq\label{kanari} \int D\left(
\sum_{i=1}^n \delta_{x_i}, \eta\right) \mu_{L,R} (d\eta) &=&
\lim_{t\to\infty}\int\E_{\eta} \caD\left( \sum_{i=1}^n \delta_{x_i},
\eta_t\right) \nu_\lambda^m (d\eta)
\nonumber\\
&=&
\lim_{t\to\infty}\int\hat{\E}^{\sip,abs}_{x_1,\ldots,x_n} \left(\caD (\sum_{i=1}^n \delta_{X_i(t)}, \eta)\right)\nu_\lambda^m (d\eta)
\nonumber\\
&\geq & \lim_{t\to\infty}\hat{\E}^{IRW(m), abs}_{x_1,\ldots,x_n} \left(\int\caD (\sum_{i=1}^n \delta_{X_i(t)}, \eta)\nu_\lambda^m (d\eta)\right)
\nonumber\\
&=&\prod_{i=1}^n \lim_{t\to\infty}\hat{\E}^{IRW(m), abs}_{x_i} \rho (X_i (t))
\nonumber\\
&=&
\prod_{i=1}^n \int D\left( \delta_{x_i}, \eta\right) \mu_{L,R} (d\eta)
\eeq
where we denoted $\hat{\E}^{\sip,abs}$ for expectation over $\sip$ particles
absorbed at $ \{ 0, N+1\}$, and $\hat{\E}^{IRW(m), abs}$ for expectation
over a system of independent random walkers (jumping at rate $2m$)
absorbed (at rate 1) at $\{ 0, N+1\}$.
\epr
\br

\ben
\item Proposition \ref{kraai} is in agreement with the findings of
\cite{gkr}, where the covariance of $\eta_i, \eta_j$ in
the measure $\mu_{L,R}$
was computed explicitly, and turned out to be
positive.
\item For the nearest neighbor SEP on $\{1,\ldots,N\}$
driven at the boundaries, we have self-duality
with absorption of dual particles
at $\{ 0,N+1\}$ and
duality function
\[
\caD_{SEP} \left(\sum_{i=1}^n\delta_{x_i}, \eta\right) =
\prod_{i=1}^n \eta_{x_i}
\]
where $\eta_0:= \rho_L, \eta_{N+1} =\rho_R$. Since for SEP
particles we have the comparison inequality of Liggett,
we have as an analogue of \eqref{noneqcor} in the SEP context,
\[
 \int \prod_{i=1}^n \eta_{x_i}\ \mu_{L,R} (d\eta) \leq \prod_{i=1}^n\int  \eta_{x_i}\ \mu_{L,R} (d\eta)
\]
i.e., $\eta_{x_i}$ are negatively correlated. The same
holds for the non-equilibrium $SEP(n)$ driven  by appropriate boundary generators.
This is in agreement with the results in \cite{spohn}, where the two-point function of the measure $\mu_{L,R}$
is computed, and with the work of \cite{der}, where some
multiple correlations are explicitly computed.
\item We expect the KMP-model, a model of heat conduction introduced and
studied in \cite{kmp} to also have positive correlations.
Indeed, the KMP and the $BEP(2)$ model are related by a so-called
instantaneous thermalization limit \cite{gkrv}. Therefore, it is natural
to think that similar correlation inequalities should hold for the KMP
as we have derived for the BEP. The limit to obtain the KMP from the BEP
is however difficult to perform on the level of the $n$-particle
representation and it is thus not clear (to us) how to prove that the KMP preserves the positive
correlation structure of the BEP. A positive hint in this direction comes
from the explicit expression of the two point function which
has been computed for the KMP in the non-equilibrium context in \cite{Bertini}.
\een
\er

\end{document}